\documentclass[a4paper,12pt]{article}
\usepackage{amsmath,fancyhdr, amssymb, amsthm, geometry, enumerate, graphicx, amsfonts, hyperref, mathrsfs}
\usepackage[latin1]{inputenc}
\usepackage{tikz}
\usetikzlibrary{shapes,arrows}
\usetikzlibrary{shadows}
\usepackage{hyperref}

\author{Neelam$^a$, \  Shah Jahan$^{a,\dag}$}
\title{ Multi-transitivity with respect to a vector in non-autonomous discrete dynamical systems }
\theoremstyle{definition}
\newtheorem{defn}{Definition}[section]
\providecommand{\keywords}[1]{\textbf{Keywords :} #1}
\providecommand{\msc}[1]{\textbf{Mathematics Subject Classification(2010)} #1}
\theoremstyle{plain}
\newtheorem{thm}{Theorem}[section]

\theoremstyle{definition}
\newtheorem{prb}{Problem}[section]
\newtheorem{exm}{Example}[section]
\newtheorem{rmk}{Remark}[section]

\begin{document}
\date{}
\maketitle

\begin{abstract} The aim of this research is to introduce the notion of multi-transitivity in non-autonomous discrete dynamical systems (NDDS) with respect to a vector. Necessary and sufficient conditions are obtained under which a NDDS is multi-transitive and strongly multi-transitive. Examples and counter-examples are constructed to justify the results in Remark 3.7 \cite{7} what is true in autonomous discrete dynamical systems (ADDS) need not be true in NDDS. Also, an affirmative answer to the question posed by Salman and Das \cite{2} is given. Finally, the Furstenberg family generated by a vector characterises multi-transitivity with respect to a vector for NDDS is obtained.
\end{abstract}


\msc{37B20; 37B55; 54B10; 54H20}
\bigskip\renewcommand{\thefootnote}{\fnsymbol{footnote}}
\footnotetext{\hspace*{-5mm}
\renewcommand{\arraystretch}{1}
\begin{tabular}{@{}r@{}p{11cm}@{}}
$^\dag$& the corresponding author. \emph{Email addresses}: shahjahan@cuh.ac.in (S. Jahan)\\
$^a$&Department of Mathematics, Central University of Haryana, Mohindergarh-123031, India
\end{tabular}}
\vspace{-2mm}
\section{Introduction}
Dynamical systems are one of the most important and useful branches of mathematics. It describes how a point in space(topological space or measurable space, or differential manifold) changes over time. Most dynamical systems are modelled using autonomous systems, in which the governing rule is considered to be constant over time. But many real-life problems, like population growth, weather forecasting, chemical reactors, the ageing of machinery, and fluid momentum, change with time. These dynamical systems can be modelled using non-autonomous discrete systems that involve a time-variant governing rule. As a result, external variables influence non-autonomous systems, which are widely used in biology, physics, engineering, and other domains for the last decades (cf. Kloeden P, cf.\cite{KloedenP, ZhangG}. In \cite{MR1402417} Kolyada and Snoha introduced and studied the concept of NDDS.\\
A NDDS is a pair $(\mathbb{D}, f_{1, \infty})$, where $f_{1, \infty}=\{f_n\}_{n=1}^\infty$ be a sequence consisting of continuous self-maps from $\mathbb{D}$ to $\mathbb{D}.$ 

In this paper,  given a topological dynamical system $(\mathbb{D}, f)$, where $\mathbb{D}$ is a compact metric space with a metric $d$ and $ f : \mathbb{D} \rightarrow \mathbb{D}$ is a continuous map, we consider the following NDDS as an analogue of non-autonomous difference equation given by:
\begin{equation}\label{1}
x_{n+1} = f_n(x_{n}), \mbox{for every natural number n,}
\end{equation}
where $f_n: \mathbb{D} \rightarrow \mathbb{D}$ is a continuous map. Note that, when $f_n = f $ then (\ref{1}) becomes the ADDS.\\
Recall that the composition $f^{n}_1$ = ${f_{n} \circ \ldots \circ f_2 \circ f_1}$ and the $n$th-iterate by $(f_{1, \infty})^{n} = \{f_{n(i-1)+1}^{n}\}_{i=1}^{\infty}$ for all natural number $n$ and $f^0_1 = Id,$ represent the identity map.\\
The orbit of a point $y \in \mathbb{D}$ is  denoted by $\mathscr{O}_{f_{1, \infty}}(y)$ and is the set
$$\mathscr{O}_{f_{1, \infty}}(y)= \{ y, f_1(y), f_1^{2}(y), \ldots, f_1^{n}(y), \ldots  \}.$$
In \cite{5} Moothathu  introduced the concept of multi-transitivity for ADDS. In \cite{7} ``Chen et.al have generalized the concept of multi-transitivity by introducing the notion of multi-transitivity with respect to(w.r.t) a vector for ADDS". They demonstrated that multi-transitivity may be defined for ADDS by hitting time sets of open sets. Recently in \cite{2} ``Salman and Das introduced and studied the notion of multi-transitivity, thick transitivity and weak mixing in NDDS". \\

In this paper, we introduce the notion of multi-transitivity in NDDS with respect to a vector. Let $(\mathbb{D}, f_{1, \infty})$ be a NDDS, it is multi-transitive if the system $(\mathbb{D}^n, (f_{1, \infty}) \times (f_{1, \infty})^{2} \times \ldots \times (f_{1, \infty})^{n})$ is transitive for every natural number $n$. We give the necessary and sufficient conditions under which a NDDS is multi-transitive and strongly multi-transitive. Examples and counter-examples are constructed to justify the results in \cite{7} what is true in ADDS need not be true in NDDS. An affirmative answer to the question posed in \cite{2} is given. Further, the Furstenberg family generated by a vector characterises multi-transitive NDDS w.r.t. a vector.\\

The outline of the paper is as follows: in Section 2, basic definitions and notations needed for following sections are given. We introduced the notion of multi-transitivity in NDDS w.r.t. a vector, and a counter-example is provided to justify that certain results related to this concept that are true in ADDS need not be true in NDDS in Section 3. Further, a positive answer to a question posed by Salman and Das is given. In Section 4, we have characterized multi-transitive NDDS w.r.t. a vector by the Furstenberg family generated by a vector.

\section{Basic Notations and Definitions}
This section includes notations and basic definitions of dynamical systems, which we will use in the subsequent sections.
Throughout this paper, we will denote $\mathbb{N}= \{1,2,3,...\},$, $\mathbb{Z}= \{...,-2,-1,0,1,2,...\},$ $\mathbb{Z}_+= \{0,1,2,...\}$,  $\mathbb{S}$ will be the unit circle and  $\{f_n\}_{n=1}^\infty =f_{1, \infty}$. We will use $\mathcal{F}$ to denote the Furstenberg family.
For any positive integers $k, m, n \in \mathbb{N}$:
\begin{center}
$(f_k)^{m}$ = ${\underbrace{f_k \circ f_k \circ \ldots \circ f_k}_\text{m-times}} \  \ , \  \ (f_k)^{-m} = {\underbrace{f_k^{-1} \circ f_k^{-1} \circ \ldots \circ f_k^{-1}}_\text{m-times}}$
\end{center}
For $p$, $n_{i}$ (i = 1, 2, \ldots , $p$) $\in \mathbb{N}$, the cartesian product $\mathbb{N}^p = {\underbrace{\mathbb{N} \times \mathbb{N} \times \cdots \times \mathbb{N}}_\text{p-times}}$ is the collection of all ordered $p$-tuples ($n_1$, $n_2$,\ldots , $n_p)$.\\
$(\mathbb{D}^n, (f_{1, \infty})^{(n)})$ defined the $n$-th product system of $(\mathbb{D}, f_{1, \infty})$, where $\mathbb{D}^n = {\underbrace{\mathbb{D} \times \mathbb{D} \times \cdots \times \mathbb{D}}_\text{n-times}}$ and $ (f_{1, \infty})^{(n)} = {\underbrace{(f_{1, \infty}) \times (f_{1, \infty}) \times \cdots \times (f_{1, \infty})}_\text{n-times}}.$\\
 Assume that $(\mathbb{D}, f_{1, \infty})$ is an NDDS. Let $A$, $B$ $\subset$ $\mathbb{D}$, the hitting time set of $A$ and $B$ is defined by $$ N(A, B) = \{n \in \mathbb{N} : f_1^{n}(A) \cap B \neq \phi\} = \{n \in \mathbb{N} : A \cap (f_1^{n})^{-1}(B) \neq \phi \}. $$
In this paper, unless otherwise stated we shall assume $A$, $B$, $A_i$ and $B_i$ for any natural number $i$ to be non-empty open subsets of $\mathbb{D}$.
\begin{defn}
Let $(\mathbb{D}, f_{1, \infty})$ be a NDDS, it is topologically transitive if for every two  $A$, $B$ $\subset$ $\mathbb{D}$, the hitting time set $ N(A, B) \neq \phi $.
\end{defn}

\begin{defn}
Let $(\mathbb{D}, f_{1, \infty})$ be a NDDS, it is totally transitive if $(\mathbb{D}, (f_{1, \infty})^{n})$ is transitive for any $n \in \mathbb{N} $.
\end{defn}

\begin{defn}
Let $(\mathbb{D}, f_{1, \infty})$ be a NDDS, it is weakly mixing if $(\mathbb{D} \times \mathbb{D}, (f_{1, \infty}) \times (f_{1, \infty}))$ is transitive.
\end{defn}

\begin{defn}
Let $(\mathbb{D}, f_{1, \infty})$ be a NDDS, it is weakly mixing of order $k$ ($k \geq$ 2) if the product system $(\mathbb{D}^k, \underbrace{(f_{1, \infty}) \times (f_{1, \infty}) \times \cdots \times (f_{1, \infty})}_\text{k - times})$ is transitive.
\end{defn}

\begin{defn}
Let $(\mathbb{D}, f_{1, \infty})$ be a NDDS, it is strongly mixing if the hitting time set $N(A, B)$ is cofinite for every two $A$, $B$ $\subset$ $\mathbb{D}$, i.e., there exists $M \in \mathbb{N}$ such that $$\{M, M + 1, \ldots\} \subset N(A, B).$$
\end{defn}

\begin{rmk}
By definitions, Strongly Mixing $\Rightarrow$ Weakly mixing of all orders $\Rightarrow$ Weakly mixing.
\end{rmk}

Let $(\mathbb{D}, f_{1, \infty})$ be a NDDS, it is multi-transitive if the system $(\mathbb{D}^n, (f_{1, \infty}) \times (f_{1, \infty})^{2} \times \ldots \times (f_{1, \infty})^{n})$ is transitive for every natural number $n$.

\begin{rmk}
By definitions, Multi-transitivity $\Rightarrow$ Total Transitivity.
\end{rmk}

Let $(\mathbb{D}, f_{1, \infty})$ be a NDDS, $(\mathbb{D}, f_{1, \infty})$  has dense small periodic sets if for every $A$ $\subset$ $\mathbb{D}$ there exists a closed subset $E$ of $A$ and for any natural number $k$ such that $(f_1^n)^{k}(E) \subset E$ i.e., $E$ is invariant for $(f_1^n)^{k}$.
If $(\mathbb{D}, f_{1, \infty})$ has dense periodic points, then it is obvious that it has dense small periodic sets.\\
``In \cite{WHX}, Huang and Ye showed
that a system which is totally transitive and has dense small periodic sets is disjoint from every minimal system." We name such a system HY-system, in shorthand.
\begin{defn}
Let $(\mathbb{D}, f_{1, \infty})$ be a NDDS,  if it is totally transitive and possesses dense small periodic sets, it is HY-system.
\end{defn}
\begin{defn}
Dynamical systems $(\mathbb{D}, f_{1, \infty})$ and $(\mathbb{Y}, g_{1, \infty})$ are called weakly disjoint if $(\mathbb{D} \times \mathbb{Y}, (f_{1, \infty}) \times (g_{1, \infty}))$ is transitive.
\end{defn}
\begin{defn}
 If an NDDS $(\mathbb{D}, f_{1, \infty})$ is weakly disjoint from any transitive system, it is called mildly mixing.
\end{defn}
\begin{defn}
Let $\mathcal{P}$ denote the power set of $\mathbb{N}$. If $\mathcal{F} \subseteq \mathcal{P}$ is is hereditary upward, it is called as \textit{Furstenberg family}, i.e
\begin{center}
$F_1 \subset F_2$ \ \  and  \ \ $F_1 \in \mathcal{F}$ \  \ $\Rightarrow $ \ \ $F_2 \in \mathcal{F}$.
\end{center}
\end{defn}
\begin{defn}
If a family $\mathcal{F}$  is a non-empty proper subset of $\mathcal{P}$, i.e., $\mathcal{F} \neq \phi$ or $\mathcal{P}$, it is called proper.
\end{defn}

\section{Main results}
We start this section by defining multi-transitivity w.r.t. a vector. Certain necessary and sufficient conditions are obtain under which $(\mathbb{D}, (f_{1, \infty})^n)$ is multi-transitive and strongly multi-transitive. Firstly, we introduce multi-transitivity in NDDS w.r.t a vector.
\begin{defn}
Let $(\mathbb{D}, f_{1, \infty})$ be a NDDS. Let a vector $\textbf{a}$ = $(a_1, a_2, \ldots, a_p)$ be  in  $\mathbb{N}^p$. We say that
\begin{enumerate}
\item if the system $(\mathbb{D}^p, (f_{1, \infty})^{(\textbf{a})})$ is transitive, then $(\mathbb{D}, f_{1, \infty})$ multi-transitive w.r.t the vector $\textbf{a}$ (or in short $\textbf{a}$-transitive), where $(f_{1, \infty})^{(\textbf{a})} = (f_{1, \infty})^{a_1} \times (f_{1, \infty})^{a_2} \times \cdots \times (f_{1, \infty})^{a_p}$;
\item if for any natural number $n$, $(\mathbb{D}, f_{1, \infty})$ is multi-transitive w.r.t (1, 2,...,$n$), then $(\mathbb{D}, f_{1, \infty})$ is multi-transitive.
\item  if for any natural number $n$, $(\mathbb{D}, f_{1, \infty})$  is multi-transitive w.r.t any vector in $\mathbb{N}^n,$ then $(\mathbb{D}, f_{1, \infty})$ is strongly multi-transitive.
\end{enumerate}
\end{defn}
\vspace{0.5cm}
\begin{thm}\label{T1}
Let $(\mathbb{D}, f_{1, \infty})$ be a NDDS. Then $(\mathbb{D}, f_{1, \infty})$ is multi-transitive iff $(\mathbb{D}, (f_{1, \infty})^n)$ is multi-transitive, for any natural number n.
\end{thm}
\begin{proof} First suppose that $(\mathbb{D}, f_{1, \infty})$ is multi-transitive. Then for any natural number $k$ and $\textbf{a}$ = (1, 2, $\ldots$, $k$), $(\mathbb{D}, f_{1, \infty})$ is $\textbf{a}$-transitive, i.e., $(\mathbb{D}^k, (f_{1, \infty})^{(\textbf{a})})$ is transitive.\\
\underline{To show:} \hspace{0.2cm} $(\mathbb{D}, (f_{1, \infty})^n)$ is $\textbf{a}$-transitive, i.e., $(\mathbb{D}^k, ((f_{1, \infty})^n)^{(\textbf{a})})$ is transitive.\\
Let $A_1$, $A_2$, $\ldots$, $A_k$, $B_1$, $B_2$, $\ldots$, $B_k$ $\subset$ $\mathbb{D}$.\\
As $(\mathbb{D}, f_{1, \infty})$ is multi-transitive. \\
Let $\textbf{a}^{'}$ = (1, 2, $\ldots$, $kn$) and $A_1^{'}$, $A_2^{'}$, $\ldots$, $A_{kn}^{'}$, $B_1^{'}$, $B_2^{'}$, $\ldots$, $B_{kn}^{'}$ be non-empty open subsets of $\mathbb{D}^{kn}$.\\
$\Rightarrow \hspace{0.5cm}(\mathbb{D}, f_{1, \infty})$ is $\textbf{a}^{'}$-transitive, i.e., $(\mathbb{D}^{kn}, (f_{1, \infty})^{(\textbf{a}^{'})})$ is transitive.\\
$$\Rightarrow \hspace{0.4cm}N_{(f_{1, \infty})^{(\textbf{a}^{'})}}(A_1^{'} \times A_2^{'} \times \ldots \times A_{kn}^{'}, B_1^{'} \times B_2^{'} \times \ldots \times B_{kn}^{'}) \neq \phi.$$ \\
Take $i$ = 1, 2, $\ldots$, $n$,\\
$$ A_{i}^{'} = A_1 \hspace{0.4cm} A_{n+i}^{'} = A_2 \hspace{0.4cm} A_{(k-1)n+i}^{'} = A_k$$ $$B_{i}^{'} = B_1 \hspace{0.4cm} B_{n+i}^{'} = B_2 \hspace{0.4cm} B_{(k-1)n+i}^{'} = B_k$$
Also,\\
$$N_{(f_{1, \infty})^{(\textbf{a}^{'})}}(A_1^{'} \times A_2^{'} \times \ldots \times A_{kn}^{'}, B_1^{'} \times B_2^{'} \times \ldots \times B_{kn}^{'}) \subset N_{((f_{1, \infty})^{n})^{(\textbf{a})}}(A_1 \times A_2 \times \ldots \times A_{k}, B_1 \times B_2 \times \ldots \times B_{k}).$$
$$\Rightarrow \hspace{0.4cm} N_{((f_{1, \infty})^{n})^{(\textbf{a})}}(A_1 \times A_2 \times \ldots \times A_{k}, B_1 \times B_2 \times \ldots \times B_{k}) \neq \phi.$$
$\Rightarrow$ \hspace{0.4cm} \hspace{0.2cm} $(\mathbb{D}^k, ((f_{1, \infty})^n)^{(\textbf{a})})$ is transitive, i.e., $(\mathbb{D}, (f_{1, \infty})^n)$ is $\textbf{a}$-transitive.\\
$\Rightarrow$ \hspace{0.4cm} \hspace{0.2cm}$(\mathbb{D}, (f_{1, \infty})^n)$ is multi-transitive.\\
\hspace{0.2cm} \textit{Sufficient:-} \hspace{0.2cm} For any natural number $k$ and $\bold{a}$ = (1, 2, $\ldots$, $k$), $(\mathbb{D}, (f_{1, \infty})^n)$ is $\textbf{a}$-transitive, i.e., $(\mathbb{D}^k, ((f_{1, \infty})^n)^{(\textbf{a})})$ is transitive.\\
\underline{To show:} \hspace{0.2cm} $(\mathbb{D}, f_{1, \infty})$ is $\textbf{a}$-transitive, i.e., $(\mathbb{D}^k, (f_{1, \infty})^{(\textbf{a})})$ is transitive.\\
As, $((f_{1, \infty})^{(\textbf{a})})^{n} = ((f_{1, \infty})^{n})^{(\textbf{a})}$ and $(\mathbb{D}^k, ((f_{1, \infty})^n)^{(\textbf{a})})$ is transitive.\\
$\Rightarrow$ $(\mathbb{D}^k, ((f_{1, \infty})^{(\textbf{a})})^{n})$ is transitive.\\
Taking $n$ = 1, $(\mathbb{D}^k, (f_{1, \infty})^{(\textbf{a})})$ is transitive, i.e., $(\mathbb{D}, f_{1, \infty})$ is $\textbf{a}$-transitive.\\
Hence, $(\mathbb{D}, f_{1, \infty})$ is multi-transitive.
\end{proof}

\begin{thm}\label{T2}
Let $(\mathbb{D}, f_{1, \infty})$ be a NDDS. Then $(\mathbb{D}, f_{1, \infty})$ is strongly multi-transitive iff $(\mathbb{D}, (f_{1, \infty})^n)$ is strongly multi-transitive, for any natural number n.
\end{thm}
\begin{proof} Suppose that $(\mathbb{D}, f_{1, \infty})$ is strongly multi-transitive.\\
For any natural number $k$ and $\textbf{a}^{'}$ = $(a^{'}_1, a^{'}_2, \ldots, a^{'}_k) \in \mathbb{N}^k$.\\
\underline{To show:} \hspace{0.2cm} $(\mathbb{D}, (f_{1, \infty})^n)$ is $\textbf{a}'$-transitive, i.e., $(\mathbb{D}^k, ((f_{1, \infty})^n)^{(\textbf{a}')})$ is transitive.\\
Since, $(\mathbb{D}, f_{1, \infty})$ is strongly multi-transitive,\\
Therefore, $(\mathbb{D}, f_{1, \infty})$ is multi-transitive w.r.t the vector $(na^{'}_1, na^{'}_2, \ldots, na^{'}_k)$.\\
$\Rightarrow$ \hspace{0.2cm}$(\mathbb{D}^k, ((f_{1, \infty})^{n})^{(\textbf{a}^{'})})$ is transitive, i.e., $(\mathbb{D}, (f_{1, \infty})^n)$ is $\textbf{a}'$-transitive.\\
$\Rightarrow$ \hspace{0.2cm}$(\mathbb{D}, (f_{1, \infty})^{n})$  is multi-transitive w.r.t the vector $\textbf{a}^{'}$.\\
Hence, $(\mathbb{D}, (f_{1, \infty})^{n})$ is strongly multi-transitive.\\

\hspace{0.2cm} \textit{Sufficient:-} \hspace{0.2cm}Suppose $(\mathbb{D}, (f_{1, \infty})^n)$ is strongly multi-transitive.\\
  For any natural number $k$ and $\textbf{a}$ = $(a_1, a_2, \ldots, a_k)$, $(\mathbb{D}, (f_{1, \infty})^n)$ is $\textbf{a}$-transitive, i.e., $(\mathbb{D}^k, ((f_{1, \infty})^n)^{(\textbf{a})})$ is transitive.\\
\underline{To show:} \hspace{0.2cm} $(\mathbb{D}, f_{1, \infty})$ is $\textbf{a}$-transitive, i.e., $(\mathbb{D}^k, (f_{1, \infty})^{(\textbf{a})})$ is transitive.\\
As, $((f_{1, \infty})^{(\textbf{a})})^{n} = ((f_{1, \infty})^{n})^{(\textbf{a})}$ and $(\mathbb{D}^k, ((f_{1, \infty})^n)^{(\textbf{a})})$ is transitive.\\
$\Rightarrow$ $(\mathbb{D}^k, ((f_{1, \infty})^{(\textbf{a})})^{n})$ is transitive.\\
Taking $n$ = 1, $(\mathbb{D}^k, (f_{1, \infty})^{(\textbf{a})})$ is transitive, i.e., $(\mathbb{D}, f_{1, \infty})$ is $\textbf{a}$-transitive.\\
Hence, $(\mathbb{D}, f_{1, \infty})$ is strongly multi-transitive.
\end{proof}
Next, following is an example which shows that $(\mathbb{D}, f_{1, \infty})$ is strongly multi-transitive.
\begin{exm}
Let $(\mathbb{D}, f_{1, \infty})$ be the NDDS, where $\mathbb{D} = \Sigma = \{ (\ldots, x_{-2}, x_{-1}, x_0, x_1, x_2, \ldots) : x_i =$ 0 or 1\} $\forall$ $i \in \mathbb{Z}$ and  $$d(x,t) = \sum_{i=-\infty}^{\infty} \frac{|x_i - t_i|}{2^{|i|}}$$ be the metric, for $x = (\ldots, x_{-2}, x_{-1}, x_0, x_1, x_2, \ldots)$, $t = (\ldots, t_{-2}, t_{-1}, t_0, t_1, t_2, \ldots) \in  \Sigma$. \\

Shift map $\sigma : \Sigma \rightarrow \Sigma$ is defined by $$\sigma((\ldots, x_{-2}, x_{-1}, \boxed{x_0}, x_1, x_2, \ldots)) = (\ldots, x_{-2}, x_{-1}, x_0, \boxed{x_1}, x_2, \ldots)$$ and $$\sigma^{-1}((\ldots, x_{-2}, x_{-1}, \boxed{x_0}, x_1, x_2, \ldots)) = (\ldots, x_{-2}, \boxed{x_{-1}} , x_0, x_1, x_2, \ldots).$$

Now taking, $f_1 = \sigma$, $f_2 = \sigma^{-4}$, $f_3 = \sigma^{4}$, so as $f_{1, \infty} = \{\sigma, \sigma^{-4}, \sigma^{4}, \sigma, \sigma^{-4}, \sigma^{4}, \ldots \}$\\

Since, $f_1^{3m}(\sigma) = \sigma^{m}$ and $\sigma$ are weakly mixing.\\
Thus, $(\mathbb{D}, f_{1, \infty})$ is also weakly mixing. Also, it is totally transitive. Therefore, $(\mathbb{D}, f_{1, \infty})$ is strongly multi-transitive.\\
By Theorem \ref{T1} $(\mathbb{D}, (f_{1, \infty})^n)$ is strongly multi-transitive.\\
By definitions, Strongly multi-transitive $\Rightarrow$ Multi-transitive.\\
$\Rightarrow \hspace{0.4cm} (\mathbb{D}, f_{1, \infty})$ is multi-transitive.\\
Hence, by Theorem \ref{T2} $(\mathbb{D}, (f_{1, \infty})^n)$ is multi-transitive.
\end{exm}

\begin{rmk}\label{R1}
\normalfont{In \cite{7} the authors have shown that for ADDS $(\mathbb{D}, f)$, $n \in \mathbb{N}$ and $\textbf{a}$ = $(a_1, a_2, \ldots, a_p)$ $\in$ $\mathbb{N}^p$. Then the subsequent statements are equivalent:
\begin{enumerate}
\item $(\mathbb{D}^k, f \times f^{2} \times \ldots \times f^{k})$ is weakly mixing, for every natural number $k$, ;
\item $(\mathbb{D}, f)$ is multi-transitive and weakly mixing.
\end{enumerate}}
\end{rmk}
\begin{rmk}\label{R2}
\normalfont{``In \cite{7} the authors have shown that for ADDS $(\mathbb{D}, f)$, if $(\mathbb{D}, f)$ is an HY-system, then $(\mathbb{D}, f)$ is strongly multi-transitive". In (\cite{7}, Remark 3.7.) author also stated that this can be proved by only showing $(\mathbb{D}, f)$ multi-transitive as every HY-system is weakly mixing for ADDS (see, \cite{WHX}, \cite{OP})}.
\end{rmk}
We provide an example to show that conditions in Remark \ref{R1} are not equivalent for NDDS. Also, Remark \ref{R2} is not true for NDDS.
\begin{exm}
\normalfont{Let $\mathbb{S}^1$ be the circle of unit radius and $(\mathbb{S}^1, f_{1, \infty})$ be the NDDS, where $f_n : \mathbb{S}^1 \rightarrow \mathbb{S}^1$ is defined as $f_{2n-1}(\theta) = p^n\theta$ and $f_{2n}(\theta) = \frac{1}{p^n}\theta$.\\
Also, $ f_{1, \infty} = \{ p\theta, \  \frac{1}{p}\theta, \  p^2\theta, \  \frac{1}{p^2}\theta, \ldots,  \ p^n\theta, \  \frac{1}{p^n}\theta, \ldots \}$}, for any $p \geq 2,$ and for every $n \in \mathbb{N}$.\\
Since $f = f_1 = p\theta$ is strongly mixing as the hitting time set $N(A, B)$ is cofinite for every two non-empty open sets $A$, $B$ $\subset$ $\mathbb{S}^1$ and $f_1^{2k-1}(\theta) = f^k(\theta)$ for every $k \in \mathbb{N}.$\\
$\Rightarrow \hspace{0.5cm} (\mathbb{S}^k, (f_{1, \infty}) \times (f_{1, \infty})^{2} \times \ldots \times (f_{1, \infty})^{k})$ is weakly mixing.\\
Now, $f_1^{2m}(\theta) = \theta$ for all $\theta \in \mathbb{S}^1$ and for all natural number $m$.\\
$\Rightarrow \hspace{0.5cm}$ For every disjoint non-empty open subsets $A, B$ of $\mathbb{S}^1$, we have $f_1^{2m}(A) \cap B = \phi$ for every $m \in \mathbb{N}$.\\
$\Rightarrow \hspace{0.5cm}$ $(\mathbb{S}^1, f_{1, \infty})$ can not be multi-transitive.
\end{exm}
From \textbf{Example 3.2.} we have seen that these conditions(\textbf{Remark 3.1.}) are not equivalent and $(\mathbb{S}^1, f_{1, \infty})$ as defined in \textbf{Example 3.2.} is not multi-transitive. Further, it is also not strongly multi-transitive.  So, \textbf{Remark 3.2.} is not true for NDDS.
\vspace{0.5cm}

Recently, Salman and Das \cite{2} obtained a sufficient condition under which a ``minimal NDDS, multi-transitivity, thick transitivity and weakly mixing of all orders are equivalent". They posed the following question:
\begin{prb}
``For an autonomous discrete dynamical system $(\mathbb{D}, f)$, it is proved that if $(\mathbb{D}, f)$ is mildly mixing, then it is multi-transitive. Does the similar result hold for a non-autonomous discrete dynamical system?"
\end{prb}
Here, we have given an affirmative answe to the problem for NDDS.
\vspace{0.5cm}
\begin{thm}\label{T4}
Let $(\mathbb{D}, f_{1, \infty})$ be a NDDS. If $(\mathbb{D}, f_{1, \infty})$ is mildly mixing, then for every $\textbf{a}$ = $(a_1, a_2, \ldots, a_p)$ $\in$ $\mathbb{N}^p$, the system $(\mathbb{D}^p, (f_{1, \infty})^{(\textbf{a})})$ is also mildly mixing. In particular, $(\mathbb{D}, f_{1, \infty})$ is strongly multi-transitive.
\end{thm}
\begin{proof} We will prove this result by the principle of  mathematical induction on length $p$ of $\textbf{a}$.\\
When $p$ = 1, $\textbf{a}$ = $(a_1=k) \in \mathbb{N}$, $(\mathbb{D}, f_{1, \infty})$ is mildly mixing, we have to show that $(\mathbb{D}, (f_{1, \infty})^{(\textbf{a})})$ i.e., $(\mathbb{D}, (f_{1, \infty})^{k})$ is also mildly mixing.\\
By fixing a natural number $k$ and let $(\mathbb{Y}, g)$ be a transitive system.\\
 We need to show $(\mathbb{D}, (f_{1, \infty})^{k})$ is also mildly mixing, i.e., $(\mathbb{D} \times \mathbb{Y}, (f_{1, \infty})^{k} \times (g_{1, \infty}))$ is transitive.\\
Let $A_1, A_2$ $\subset$ $\mathbb{D}$ and $B_1, B_2$ be non-empty open subsets of $\mathbb{Y}$.\\
Let $y$ is a transitive point of $(\mathbb{Y}, g_{1, \infty})$.\\
Let $\mathcal{K}$ = \{0, 1, $\ldots$, $k$ - 1\} be a system with discrete topology and $\mathcal{Z} = \mathbb{Y} \times \mathcal{K}$.\\
We define $h_{m} \ : \ \mathcal{Z} \rightarrow \mathcal{Z}$ by
\begin{center}
$h_{m}(x, i)$ = $ \left\{ \begin{array}{ll} (x, i + 1), & i = 0, \ldots, k - 2, \\ (g_{m}(x), 0), & i = k - 1.
\end{array} \right. $
\end{center}
Taking $x = y$ and $i = 0$, we have
\begin{center}
$h_{m}(y, 0)$ = $(y, 1)$
\end{center}
\begin{center}
	$h_{m}(y,1)$ = $(y, 2)$
\end{center}
\begin{center}
$\vdots$
\end{center}
\begin{center}
	$h_{m}(y, k - 2)$ = $(y, k - 1)$
\end{center}
\begin{center}
	$h_{m}(y, k - 1)$ = $(g_m(y), 0)$
\end{center}
Thus, we get $\overline{\mathscr{O}(y, 0)}$ = $\mathcal{Z}$.\\
$\Rightarrow$ $(\mathcal{Z}, h_{1, \infty})$ is a transitive system with $(y, 0)$ as a transitive point.

\noindent Take $B_1^{'}$ = $B_1 \times$ \{0\} and $B_2^{'}$ = $B_2 \times$ \{0\}\\
As $(\mathbb{D}, f_{1, \infty})$ is mildly mixing and $(\mathcal{Z}, h_{1, \infty})$ is transitive,\\
$\Rightarrow$ $(\mathbb{D} \times \mathcal{Z} , (f_{1, \infty}) \times (h_{1, \infty}))$ is transitive.\\
Then $\exists$ a natural number $n$ which satisfy
\begin{center}
$A_1 \ \cap (f_1^{n})^{-1}(A_2) \ \neq \phi$ \ \ and \ \  $B_1^{'} \ \cap \ (h_1^{n})^{-1}(B_2^{'}) \neq \phi$
\end{center}
By computation of $h_{m}$, we have $k \vert n$. Take $l = n/k$.\\
Then,
\begin{center}
$A_1 \cap (f_1^{kl})^{-1}(A_2) \neq \phi$ \ \ and \ \ $B_1 \cap (g_1^{l})^{-1}(B_2) \neq \phi,$
\end{center}
$\Rightarrow$ $(\mathbb{D} \times \mathbb{Y}, (f_{1, \infty})^k \times (g_{1, \infty}))$ is transitive.\\
Now, for any natural number $p > 1,$ we suppose that the result holds for $p$ - 1. So, we have $(\mathbb{D}^{p-1}, (f_{1, \infty})^{(\textbf{a}^{'})})$ is mildly mixing where \textbf{a}$^{'}$ = $(a_1, a_2, \ldots, a_{p - 1})$.\\
Let \textbf{a} = $(a_1, a_2, \ldots, a_{p - 1}, a_p)$ be a vector in $\mathbb{N}^{p}$ with length $p$ and $(\mathbb{Y}, g_{1, \infty})$ be a transitive system.\\
 Since, $(\mathbb{D}, (f_{1, \infty})^{a_p})$ is a mildly mixing system, so, $(\mathbb{D} \times \mathbb{Y}, (f_{1, \infty})^{a_p} \times (g_{1, \infty}))$ is transitive.\\
Also, $(\mathbb{D}^{p-1}, (f_{1, \infty})^{(\textbf{a}^{'})})$ is mildly mixing system,\\
 we have $(\mathbb{D}^{p - 1} \times \mathbb{D} \times \mathbb{Y}, (f_{1, \infty})^{(\textbf{a}^{'})} \times (f_{1, \infty})^{a_p} \times (g_{1, \infty}))$ is transitive i.e., $(\mathbb{D}^p \times Y, (f_{1, \infty})^{(\textbf{a})} \times (g_{1, \infty}))$ is transitive.\\
$\Rightarrow$ \hspace{2cm} $(\mathbb{D}^p, (f_{1, \infty})^{(\textbf{a})})$ is mildly mixing.\\
Hence, the result is true for the length $p$ of \textbf{a}.
\end{proof}
We give the following example to support the Theorem \ref{T4}.
\begin{exm}
\normalfont{Let $(\mathbb{D}, f_{1, \infty})$ be the NDDS and $(\mathbb{Y}, g_{1, \infty})$ is transitive system, where $\mathbb{D} = \mathbb{Y} = \Sigma = \{ (x_0, x_1, \ldots) : x_i =$ 0 or 1\}  $\forall$ $i \in \mathbb{Z}$ with metric $$d(x,t) = \sum_{i=-\infty}^{\infty} \frac{|x_i - t_i|}{2^{|i|}}$$ for any $x = (\ldots, x_{-2}, x_{-1}, x_0, x_1, x_2, \ldots) \in \Sigma $ and $t = (\ldots, t_{-2}, t_{-1}, t_0, t_1, t_2, \ldots) \in  \Sigma.$ \\

Define the map $\sigma : \Sigma \rightarrow \Sigma$ by $$\sigma((\ldots, x_{-2}, x_{-1}, \boxed{x_0}, x_1, x_2, \ldots)) = (\ldots, x_{-2}, x_{-1}, x_0, \boxed{x_1}, x_2, \ldots)$$ and $$\sigma^{-1}((\ldots, x_{-2}, x_{-1}, \boxed{x_0}, x_1, x_2, \ldots)) = (\ldots, x_{-2}, \boxed{x_{-1}} , x_0, x_1, x_2, \ldots).$$

Now, taking $f_1 = \sigma$, $f_2 = \sigma^{-2}$, $f_3 = \sigma^{2}$, so as $f_{1, \infty} = \{\sigma, \sigma^{-2}, \sigma^{2}, \sigma, \sigma^{-2}, \sigma^{2}, \ldots \}$}, and taking $g_1 = \sigma$, $g_2 = \sigma^{-3}$, $g_3 = \sigma^{3}$, so as $g_{1, \infty} = \{\sigma, \sigma^{-3}, \sigma^{3}, \sigma, \sigma^{-3}, \sigma^{3}, \ldots \}$\\
As $f_1^{3m}(\sigma) = \sigma^m$, $g_1^{3m}(\sigma) = \sigma^m$ and $\sigma$ is weakly mixing.\\
So, for $a_p \in \mathbb{N}$, $(\mathbb{D} \times \mathbb{Y}, (f_{1, \infty})^{a_p} \times (g_{1, \infty}))$ is mildly mixing.\\
Also, $(\mathbb{D}, f_{1, \infty})$  with respect to $\textbf{a}^{'}$ = $(a_1, a_2, \ldots, a_{p-1})$ $\in$ $\mathbb{N}^{p-1}$ is transitive.\\
Thus, we have $(\mathbb{D}^{p - 1} \times \mathbb{D} \times \mathbb{Y}, (f_{1, \infty})^{(\textbf{a}^{'})} \times (f_{1, \infty})^{a_p} \times (g_{1, \infty}))$ is transitive i.e., $(\mathbb{D}^p \times Y, (f_{1, \infty})^{(\textbf{a})} \times (g_{1, \infty}))$ is transitive, where $\textbf{a}$ = $(a_1, a_2, \ldots, a_{p})$ $\in$ $\mathbb{N}^p$.\\
Hence,  $(\mathbb{D}^p, (f_{1, \infty})^{(\textbf{a})})$ is mildly mixing.
\end{exm}
\section{Characterization of multi-transitivity w.r.t a Vector}
Gottschalk and Hedhund \cite{6} described autonomous discrete dynamical systems by using families, which was further developed by Furstenberg \cite{1}. In 2018, Vasisht and Das \cite{4} introduced the notion of Furstenberg family for NDDS. Here, firstly we define Furstenberg family by the hitting time set of open sets for NDDS.
\begin{defn}
Let $(\mathbb{D}, f_{1, \infty})$ be a NDDS and $\mathcal{F}$ be a family, if for every $A, B \subset \mathbb{D}$,  $N(A, B) \in \mathcal{F},$ then $(\mathbb{D}, f_{1, \infty})$ is called $\mathcal{F}$ - transitive.
\end{defn}
\begin{defn}
Let $(\mathbb{D}, f_{1, \infty})$ be a NDDS and $\mathcal{F}$ be a family, if $(\mathbb{D} \times \mathbb{D}, (f_{1, \infty}) \times (f_{1, \infty}))$ is $\mathcal{F}$ - transitive, then $(\mathbb{D}, f_{1, \infty})$ is called $\mathcal{F}$ - mixing.
\end{defn}
\begin{defn}
Let $\bold{a} = (a_1, a_2, \ldots, a_p)$ be a vector $\in \mathbb{N}^p$. The family $\mathcal{F}[\textbf{a}]$ generated by the vector $\textbf{a}$ is defined as
\begin{center}
\{$F \ \subset \mathbb{N} :$ for every $\textbf{n} = (n_1, n_2, \ldots, n_p) \in \mathbb{Z}_{+}^p$,$\exists$ a natural number $m$ such that $m\textbf{a} + \textbf{n} \in F^p$ \},
\end{center}
where $m\textbf{a} + \textbf{n} = (ma_1 + n_1, ma_2 + n_2, \ldots, ma_p + n_p)$.
\end{defn}
\begin{thm}\label{T3}
Let $(\mathbb{D}, f_{1, \infty})$ be a NDDS, and $\textbf{a} = (a_1, a_2, \ldots, a_p) \in \mathbb{N}^p$. Then $(\mathbb{D}, f_{1, \infty})$ is $\textbf{a}$-transitive iff it is $\mathcal{F}[\textbf{a}]$-transitive.
\end{thm}
\begin{proof} First suppose that $(\mathbb{D}, f_{1, \infty})$ is $\textbf{a}$-transitive. Let $A, B \subset \mathbb{D}$.\\
\underline{To show:}\hspace{1cm} $N(A, B) \in \mathcal{F}[\textbf{a}]$.\\
Let $(n_1, n_2, \ldots, n_p) \in \mathbb{Z}_{+}^p$.\\
As $(\mathbb{D}, f_{1, \infty})$ is $\textbf{a}$-transitive, i.e., $(\mathbb{D}^p, (f_{1, \infty})^{(\textbf{a})})$ is transitive.\\
Then, $\exists$ a natural number $m$ such that
\begin{center}
$m \in N_{(f_{1, \infty})^{(\bold{a})}}(A \times A \times \cdots \times A, (f_1^{n_1})^{-1}(B) \times (f_1^{n_2})^{-1}(B) \times \cdots \times (f_1^{n_p})^{-1}(B))$\\
\end{center}
$\Rightarrow \hspace{0.5cm} A \ \cap \ (f_1^{m\textbf{a} + \textbf{n}})^{-1}(B) \neq \phi $\\
$\Rightarrow \hspace{0.5cm} m\textbf{a} + \textbf{n} \subset N(A, B)$\\
$\Rightarrow \hspace{0.5cm} \{ma_1 + n_1, ma_2 + n_2, \ldots, ma_p + n_p\} \subset N(A, B)$.\\
Hence, $N(A, B) \in \mathcal{F}[\textbf{a}]$\\

\hspace{0.2cm} \textit{Sufficient:-} \hspace{0.2cm} Suppose that $(\mathbb{D}, f_{1, \infty})$ is $\mathcal{F}[\textbf{a}]$-transitive. Let $A_1, A_2, \ldots, A_p$ and  $B_1, B_2, \ldots, B_p$ $\subset$ $\mathbb{D}$.\\
 Using transitivity of $(\mathbb{D}, f_{1, \infty})$ to $A_{p-1}$ and $A_p$, we choose a $m_{p-1} \in \mathbb{N}$ such that
\begin{center}
$A_{p-1} \ \cap \ (f_1^{n})^{-m_{p-1}}(A_p) \neq \phi$
\end{center}
Again using transitivity of $(\mathbb{D}, f_{1, \infty})$ to $A_{p-2}$ and $A_{p-1} \ \cap \ (f_1^{n})^{m_{p-1}}(A_p)$, we choose a $m_{p-2} \in \mathbb{N}$ such that
\begin{center}
$A_{p-2} \ \cap \ (f_1^{n})^{-m_{p-2}}(A_{p-1} \ \cap \ (f_1^{n})^{m_{p-1}}(A_p)) \neq \phi$
\end{center}
Continuing the process $p-1$ times, we get a sequence $\{m_i\}_{i=1}^{p-1}$ of positive integers and
$$\begin{array}{lll}
A^{'} & = & A_1 \cap (f_1^{n})^{-m_1}(A_2 \cap (f_1^{n})^{-m_2}(\cdots \cap (A_{p-1} \cap (f_1^{n})^{-m_{p-1}}(A_p))))\\
& = & A_1 \cap (f_1^{n})^{-m_1}(A_2) \cap (f_1^{n})^{-(m_1 + m_2)}(A_3) \cap \cdots \cap (f_1^{n})^{-(m_1 + m_2 + \cdots + m_{p-1})}(A_p) \neq \phi
\end{array}$$
Similarly, we can obtain a sequence of $\{m_i^{'}\}_{i = 1}^{p-1}$ with $m_i^{'} < m_i$ for $i$ = 1, 2, $\ldots$, $p-$1 and
$$\begin{array}{lll}
B^{'} & = & B_1 \cap (f_1^{n})^{-m^{'}_1}(B_2 \cap (f_1^{n})^{-m^{'}_2}(\cdots \cap (B_{p-1} \cap (f_1^{n})^{-m^{'}_{p-1}}(B_p))))\\
& = & B_1 \cap (f_1^{n})^{-m^{'}_1}(B_2) \cap (f_1^{n})^{-(m^{'}_1 + m^{'}_2)}(B_3) \cap \cdots \cap (f_1^{n})^{-(m^{'}_1 + m{'}_2 + \cdots + m{'}_{p-1})}(B_p) \neq \phi
\end{array}$$
By supposition, $N_{(f_{1, \infty})}(A^{'}, B^{'}) \in \mathcal{F}[\textbf{a}]$. Then for the vector (0, $m_1 - m^{'}_1$, $m_1 + m_2 - m^{'}_1 - m^{'}_2$, $\ldots$, $m_1 + m_2 + \cdots + m_{p-1} - m_1^{'} - m_2^{'} - \cdots - m_{p-1}^{'}$), there exist $l \in \mathbb{N}$ such that
$$\{a_1l, a_2l + m_1 -m_1^{'} \ldots, a_pl + m_1 + m_2 + \cdots + m_{p-1} - m_1^{.} - m_2^{'} - \cdots - m_{p-1}^{'} \} \subset N_{(f_{1, \infty})}(A^{'}, B^{'}).$$
Then \\
$A_1 \cap ((f_1^{n})^{a_1})^{-l}(B_1) \neq \phi$,\\
$f_1{^{n}}(A_2) \cap ((f_1^{n})^{a_2})^{-l}((f_1^{n})^{-m_1}(B_2)) \neq \phi$,\\
$\vdots$\\
${f_1^n}^{(m_1 + m_2 + \cdots + m_{p-1})}(A_p) \cap ((f_1^{n})^{a_p})^{-l}((f_1^{n})^{-(m_1 + m_2 + \cdots + m_{p-1})}(B_p)) \neq \phi$.\\
$\Rightarrow$ \hspace{0.5cm} $l \in N_{(f_{1, \infty})^{(\bold{a})}}(A_1 \times A_2 \times \cdots \times A_p, B_1 \times B_2 \times \cdots \times B_p).$\\
Hence, $(\mathbb{D}, f_{1, \infty})$ is $\textbf{a}$-transitive.
\end{proof}

Let $\mathcal{F}[\infty] = \cap_{p = 1}^{\infty} \mathcal{F}[a_p]$ , where $\bold{a}_p = (1, 2, \ldots, p)$ for $p \in \mathbb{N}$.\\
Recalling from \cite{7}, $F \subseteq \mathbb{N}$ is in $\mathcal{F}[\infty]$ iff for every natural number $p$ and every $\textbf{n} \in \mathbb{Z}_{+}^{p}$, $\exists$ $m \in \mathbb{N}$ such that $m \textbf{a}_p$ + \textbf{n} \  $\in$ $F^p$.
\vspace{0.5cm}
\begin{thm}
Let $(\mathbb{D}, f_{1, \infty})$ be a NDDS. Then $(\mathbb{D}, f_{1, \infty})$ is multi-transitive iff it is $\mathcal{F}[\infty]$-transitive.
\end{thm}
\begin{proof} For every $p \in \mathbb{N}$, $\textbf{a}_p = (1, 2, \ldots, p) \in \mathbb{N}^p$.
$(\mathbb{D}, f_{1, \infty})$ is multi-transitive, i.e., $(\mathbb{D}, f_{1, \infty})$ is $\textbf{a}_p$-transitive $\forall$ $p \in \mathbb{N}$.\\
Then, by Theorem \ref{T3}, $(\mathbb{D}, f_{1, \infty})$ is $\mathcal{F}[\textbf{a}_p]$-transitive $\forall$ $p \in \mathbb{N}$. Then for every two $A, B$ $\subset$ $\mathbb{D}$, $N(A, B) \in \mathcal{F}[\textbf{a}_p]$ $\forall$ $p \in \mathbb{N}$.\\
$\Rightarrow \hspace{0.5cm} N(A, B) \in \mathcal{F}[\infty]$.\\
$\Rightarrow \hspace{0.5cm} (\mathbb{D}, f_{1, \infty})$ is $\mathcal{F}[\infty]$-transitive.\\

  Conversely, suppose that $(\mathbb{D}, f_{1, \infty})$ is $\mathcal{F}[\infty]$-transitive, i.e., $(\mathbb{D}, f_{1, \infty})$ is $\mathcal{F}[\textbf{a}_p]$-transitive $\forall$ $p \in \mathbb{N}$.\\
Then, by Theorem \ref{T3}, $(\mathbb{D}, f_{1, \infty})$ is $\textbf{a}_p$-transitive for all natural number $p$.
Thus, $(\mathbb{D}, f_{1, \infty})$ is multi-transitive.
\end{proof}
	\textbf{Funding:} No funding.\\
\textbf{Conflict of interest:} We declare that there is no conflict of interest for
the publication of this paper.\\
\textbf{Author Contributions:} All the authors have made equal contributions in
this article. All authors read and approved the final manuscript.\\
\textbf{Availability of data and material:}  No data is used.\\
\textbf{Acknowledgements:} We would like to thank Central University of Haryana for providing necessary facilities to carry out this research and also thankful to all the persons who have made
substantial contributions to the work reported in the manuscript (e.g.,
technical help, writing and editing assistance, general support).

\end{document}